\documentclass[a4paper]{amsart}

\usepackage{amsmath,amssymb,amsthm,stmaryrd}

\usepackage{color,graphicx}
\usepackage[dvipsnames]{xcolor}
\usepackage{accents}
\usepackage{enumerate}
\usepackage{varwidth}

\usepackage[textsize=footnotesize,color=yellow!40, bordercolor=white]{todonotes}

\newtheorem{theorem}{Theorem}[section]
\newtheorem{conjecture}[theorem]{Conjecture}
\newtheorem{lemma}[theorem]{Lemma}
\newtheorem{fact}[theorem]{Fact}

\newtheorem*{claim*}{Claim}
\newtheorem*{subclaim*}{Subclaim}

\newtheorem{question}[theorem]{Question}

\newcommand{\thistheoremname}{}
\newtheorem*{genericthm}{\thistheoremname}
\newenvironment{namedthm}[1]
  {\renewcommand{\thistheoremname}{#1}%
   \begin{genericthm}}
  {\end{genericthm}}

\theoremstyle{definition}
\newtheorem{definition}[theorem]{Definition}

\newtheorem*{definition*}{Definition}
\newcounter{ca}[cl]

\newcounter{sca}[ca]

\newcounter{ssca}[sca]

\theoremstyle{remark}
\newtheorem{remark}[theorem]{Remark}
\newtheorem*{remark*}{Remark}

\usepackage{tikz}
\usetikzlibrary{decorations.pathmorphing,shapes}
\usetikzlibrary{positioning,arrows}

\usepackage[
    colorlinks=true, 
    citecolor=blue,
    linkcolor=blue,
    urlcolor=blue]{hyperref}
    
\newcommand{\ZFC}{\ensuremath{\operatorname{ZFC}} }
\newcommand{\ZF}{\ensuremath{\operatorname{ZF}} }
\newcommand{\GCH}{\ensuremath{\operatorname{GCH}} }

\newcommand{\CH}{\ensuremath{\operatorname{CH}} }
\newcommand{\OR}{\ensuremath{\operatorname{Ord}} }
\newcommand{\Ord}{\ensuremath{\operatorname{Ord}} }
\newcommand{\HOD}{\ensuremath{\operatorname{HOD}} }

\newcommand{\cf}{\ensuremath{\operatorname{cf}} }

\newcommand{\AC}{\ensuremath{\operatorname{AC}} }

\newcommand{\AD}{\ensuremath{\operatorname{AD}} }
\newcommand{\PFA}{\ensuremath{\operatorname{PFA}} }

\newcommand{\WF}{\ensuremath{\operatorname{WF}} }

\newcommand{\trcl}{\ensuremath{\operatorname{trcl}} }
\newcommand{\uB}{{\ensuremath{\sf{uB}}} }
\newcommand{\uBp}{{\ensuremath{\sf{uBp}}} }

\newcommand{\DC}{\ensuremath{\operatorname{DC}} }

\newcommand{\Col}{\ensuremath{\operatorname{Col}} }

\newcommand{\Ult}{\ensuremath{\operatorname{Ult}} }

\newcommand{\meas}{\ensuremath{\operatorname{meas}} }

\newcommand{\Def}{\ensuremath{\operatorname{Def}} }
    
\newcommand{\powerset}{{\wp}}

\newcommand{\BS}{{}^\omega\omega}

\newcommand{\PI}{\boldsymbol\Pi}
\newcommand{\SIGMA}{\boldsymbol\Sigma}

\newcommand{\bR}{\mathbb{R}}
\newcommand{\bP}{\mathbb{P}}

\newcommand{\card}[1]{{\vert #1 \vert} }

\newcommand{\Sealing}{\mathsf{Sealing}}
\newcommand{\UltimateL}{\mathsf{V = Ultimate-L}}

\newcommand{\Hom}{\ensuremath{\operatorname{Hom}}}

\newcommand{\shortpwimg}{"}

\def\k{\kappa}
\def\a{\alpha}

\begin{document}
\title{G\"odel's Program in Set Theory}



\author{Sandra M\"uller} \address{Sandra M\"uller, Institut f\"ur Diskrete Mathematik und Geometrie, TU Wien, Wiedner Hauptstra{\ss}e 8-10/104, 1040 Wien, Austria.}
\email{sandra.mueller@tuwien.ac.at}
\author{Grigor Sargsyan} \address{Grigor Sargsyan, IMPAN, Antoniego Abrahama 18, 81-825 Sopot, Poland.}
\email{gsargsyan@impan.pl}

\thanks{This research was funded in part by the Austrian Science Fund (FWF) [10.55776/I6087, 10.55776/Y1498]. For open access purposes, the authors have applied a CC BY public copyright license to any author accepted manuscript version arising from this submission. The second author's work is funded by the National Science Centre, Poland under the Maestro 15 call in the Maestro programme,
registration number UMO-2023/50/A/ST1/00258.}




\begin{abstract}
Gödel proved in the 1930s in his famous Incompleteness Theorems that not all statements in mathematics can be proven or disproven from the accepted $\ZFC$ axioms. A few years later he showed the celebrated result that Cantor's Continuum Hypothesis is consistent. Afterwards, Gödel raised the question whether, despite the fact that there is no reasonable axiomatic framework for all mathematical statements, natural statements, such as Cantor’s Continuum Hypothesis, can be decided via extending $\ZFC$ by large cardinal axioms. While
this question has been answered negatively, the problem of finding good axioms that decide natural mathematical statements remains open. There is a
compelling candidate for an axiom that could solve Gödel’s problem: $\UltimateL$. In addition, due to recent results the $\Sealing$ scenario has gained a lot of attention. We describe these candidates as well as their impact and relationship.
  \end{abstract}
\maketitle
\setcounter{tocdepth}{1}

\section{Introduction}

G\"{o}del's Program is a major set theoretic program addressing the most fundamental set theoretic issue, \textit{independence}: the inability of the basic axioms of set theory (Zermelo-Fraenkel set theory with Choice, $\ZFC$) or any (consistent) extension thereof to decide natural questions about infinite sets. Among these natural questions that cannot be decided by $\ZFC$ is Cantor's Continuum Hypothesis ($\CH)$. Deciding its truth is arguably among the most basic questions in set theory. It can informally be stated as the simple question: How many real numbers are there?

G\"odel's Program aims at removing undecidability from foundations by developing and studying natural extensions of $\ZFC$. The goal of the program is to remove incompleteness in the natural foundational theory $T$  by passing to a stronger theory that is as natural as $T$ itself and  decides all or some of the undecidable questions of $T$. G\"odel's idea outlined in \cite{Godel} was that iterating this process would resolve all undecidable questions of $\ZFC$, and that they would all be resolved within the \emph{Large Cardinal Hierarchy}. As demonstrated by the celebrated Levy-Solovay theorem, the Large Cardinal Hierarchy cannot achieve G\"{o}del's dream. At the same time, other set theoretic hierarchies, such as forcing axioms or determinacy axioms, that are deeply connected to the Large Cardinal Hierarchy, can go extremely far in deciding natural questions, including $\CH$, that are undecidable within $\ZFC$.

\section{Natural hierarchies extending the standard axioms for mathematics}

We describe the three most successful hierarchies of axioms extending $\ZF$ or $\ZFC$, the standard axioms for mathematics: the Large Cardinal Hierarchy, Determinacy Axioms, and Forcing Axioms. Determinacy assumptions, large cardinal axioms, forcing axioms, and their consequences are widely used and have many fruitful implications in set theory and even in other areas of mathematics such as algebraic topology \cite{CSS05}, topology \cite{Ny80,Fl82,CMM1}, algebra \cite{EM02}, and operator algebras \cite{Fa11}. Many such applications, in particular, proofs of consistency strength lower bounds, exploit the interplay of large cardinals and determinacy axioms. Thus, understanding the connections between determinacy assumptions, the hierarchy of large cardinals, and forcing axioms is vital to answer questions left open by $\ZFC$ itself.

\subsection*{Large Cardinal Hierarchy} Large cardinal axioms postulate the existence of large infinities. They are naturally ordered in a hierarchy given by their strength and while smaller large cardinals, such as inaccessible cardinals or Mahlo cardinals, can still exist in Gödel's constructible universe $L$, larger large cardinals such as measurable cardinals contradict the Axiom of Constructibility $V = L$.   An excellent account on large cardinals including their history and definitions can be found in \cite{Ka08}. We only outline the definitions of the three most important ones for our context here: measurable cardinals, Woodin cardinals, and supercompact cardinals.

Measurable cardinals were originally introduced by Ulam \cite{Ul30} in the context of two-valued measures. In this sense, an uncountable cardinal $\kappa$ is measurable if there is a non-trivial $\kappa$-additive, two-valued measure on $\powerset(\kappa)$. Equivalently, an uncountable cardinal $\kappa$ is measurable if there is a $\kappa$-complete, non-principal ultrafilter on $\kappa$. Measurable cardinals were further developed by Keisler and Scott \cite{Kei62, Sc61} and their equivalent characterization via elementary embeddings turns out to be very useful in set theory.

\begin{theorem}[Keisler, Scott, \cite{Kei62, Sc61}]
    An uncountable cardinal $\kappa$ is measurable if, and only if, there is a non-trivial elementary embedding \[ j \colon V \rightarrow M \] from the universe $V$ into a transitive class $M$ with critical point $\kappa$, i.e., $j(\alpha) = \alpha$ for all ordinals $\alpha < \kappa$ and $j(\kappa) > \kappa$. 
\end{theorem}

Using this result as a starting point and requiring more similarity between $V$ and $M$ results in stronger large cardinal properties. To introduce Woodin cardinals, we first define strong cardinals.

\begin{definition}
    Let $\eta$ be an ordinal. 
    \begin{enumerate}
        \item An uncountable cardinal $\kappa$ is \emph{$\eta$-strong} if there is a non-trivial elementary embedding \[ j \colon V \rightarrow M \] from the universe $V$ into a transitive class $M$ with critical point $\kappa$ such that $j(\kappa) > \eta$ and $V_{\kappa+\eta} \subseteq M$.
        \item For any set $A$, an uncountable cardinal $\kappa$ is \emph{$A$-$\eta$-strong} if there is a non-trivial elementary embedding \[ j \colon V \rightarrow M \] from the universe $V$ into a transitive class $M$ with critical point $\kappa$ such that $j(\kappa) > \eta$, $V_{\kappa+\eta} \subseteq M$, and \[ A \cap V_{\kappa+\eta} = j(A) \cap V_{\kappa+\eta}. \]
    \end{enumerate}
\end{definition}

Now we are ready to define Woodin cardinals. As the name suggests, they have been isolated by Woodin. Their definition is of a slightly different spirit as many other large cardinal properties as, for a Woodin cardinal $\delta$, it does not postulate the existence of elementary embeddings with critical point $\delta$ but instead requires that there are many large cardinals below $\delta$. In fact, a Woodin cardinal $\delta$ itself need not even be inaccessible or measurable, even though there will be many inaccessible and measurable cardinals below any Woodin cardinal. Woodin cardinals turned out to be extremely useful in set theory and have very natural deep connections to descriptive set theory, for example, in the projective hierarchy.

\begin{definition}[Woodin]
    An uncountable cardinal $\delta$ is \emph{Woodin} if for any set $A \subseteq V_\delta$, there is a $\kappa < \delta$ that is $A$-$\eta$-strong for every $\eta < \delta$.
\end{definition}

Another well-studied large cardinal notion is a supercompact cardinal. This is beyond the scope of current inner model theory in the sense that it is not known how to construct canonical inner models with supercompact cardinals. However, supercompact cardinals will be heavily used in Section \ref{sec:Sealing} of this article when we discuss $\Sealing$.

\begin{definition}
 An uncountable cardinal $\kappa$ is \emph{supercompact} if for any $\eta \geq \kappa$ there is a non-trivial elementary embedding \[ j \colon V \rightarrow M \] from the universe $V$ into a transitive class $M$ with critical point $\kappa$ such that $j(\kappa) > \eta$ and $M^{\eta} \subseteq M$. 
\end{definition}

\subsection*{Determinacy Hierarchy} Determinacy axioms postulate the existence of winning strategies in two-player games \cite{GS53, MySt62}. For every set of infinite sequences of natural numbers\footnote{As common in set theory, we will tacitly identify infinite sequences of natural numbers with reals.} $A \subseteq \mathbb{N}^{\mathbb{N}}$, we consider an infinite two-player game where player $\mathrm{I}$ and player
$\mathrm{II}$ alternate playing natural numbers $n_0, n_1, \dots$, as follows: \vspace{-0.2cm}
\[ \begin{array}{c|ccccc} \mathrm{I} & n_0 & & n_2 & &\hdots \\ \hline
    \mathrm{II} & & n_1 & & n_3 & \hdots 
   \end{array} \]

Then player $\mathrm{I}$ wins the game if and only if the sequence $x = (n_0,n_1,\dots)$ of natural numbers produced during a run of the game is
an element of $A$; otherwise, player $\mathrm{II}$ wins. We call $A$ the \emph{payoff set} of this game. A set $A$ is called \emph{determined} if and only if
one of the two players has a winning strategy in the game with payoff set $A$, meaning that there is a method by
which they can win in the game described above, no matter what their opponent does. The \textbf{Axiom of Determinacy} ($\AD$) is the statement that all sets of reals are determined. Determinacy hypotheses are known to 
enhance sets of real numbers with a great deal of canonical structure. 
But while a pioneering result of Gale and Stewart shows that every open and every closed set in $\mathbb{N}^{\mathbb{N}}$ is determined in $\ZFC$, determinacy for more complex, for example, analytic, sets cannot be proven in $\ZFC$ alone. It should be mentioned at this point that the full Axiom of Determinacy as stated above contradicts the Axiom of Choice, so it can only be considered as an extension of $\ZF$, not of $\ZFC$.

\subsection*{Forcing Axiom Hierarchy} Forcing Axioms postulate that certain sets used to build forcing extensions of the set theoretic universe exist. Forcing extensions go back to Cohen's famous proof of the consistency of the failure of the Continuum Hypothesis and are nowadays omnipresent in set theory. The first forcing axioms that was isolated was Martin's Axiom introduced by Martin and inspired by Solovay's and Tennenbaum's proof of the independence of Souslin's Hypothesis \cite{MaSo70, SoTe71}. Other well-known examples with many applications are the Proper Forcing Axiom introduced by Baumgartner (see \cite{Mo10}) and the in some sense ``ultimate forcing axiom'' Martin's Maximum introduced by Foreman, Magidor, and Shelah \cite{FMS88}. We will not go into details here as Forcing Axioms will not be in the focus of the article.

There are other well-studied hierarchies of axioms that extend $\ZFC$. Another prominent one is the hierarchy of generic embeddings and their associated ideals. Connecting this hierarchy to the three hierarchies mentioned above has led to several deep proofs as well as many interesting open questions.

\subsection*{Example: Cantor's Continuum Problem}

Arguably the most famous statement that is known to be independent from $\ZFC$ is Cantor's Continuum Problem. It was formulated by Cantor in 1878 and appeared as the first item on Hilbert's list of problems announced at the International Congress of Mathematicians in Paris in 1900. Informally, it can be phrased as the question how many real numbers there are. Or, a bit more formally, as the following question: Is there a set $A$ of size strictly between the size of the set of natural numbers $|\mathbb{N}|$ and the size of the set of real numbers $|\mathbb{R}|$? I.e., is there a set $A$ such that \[ |\mathbb{N}| < |A| < |\mathbb{R}|? \]

\begin{definition}\label{def:CH}
    For a collection of sets of reals $\Gamma$ we say that the \emph{Continuum Hypothesis holds for $\Gamma$} if and only if there is no set $A \in \Gamma$ such that $|\mathbb{N}| < |A| < |\mathbb{R}|$. In particular, the Continuum Hypothesis holds if it holds for $\powerset(\mathbb{R})$.
\end{definition}

Gödel showed in the late 1930's that there are models of set theory in which there is no such set $A$. Many years later in 1963 Cohen was able to show that there are also models of set theory in which there is such a set $A$, thereby establishing the independence of the Continuum Problem from $\ZF$. He was awarded the Fields Medal in 1966 for this outstanding result.

We are interested in the question if and how canonical extensions of $\ZF$ decide the Continuum Problem. Interestingly, the three lines of extensions of $\ZF$ described above give different answers to the Continuum Problem.

When phrased as above, determinacy gives a positive answer to the Continuum Problem:

\begin{fact}
    If the Axiom of Determinacy holds, then the perfect set property holds, i.e., every set of reals is either countable or contains a perfect subset. In particular, there is no set $A$ such that $|\mathbb{N}| < |A| < |\mathbb{R}|$.
\end{fact}

\begin{remark}
    The status of the Continuum Problem under the Axiom of Determinacy is not as clear as it might seem here. There are several formulations of $\CH$ that are all equivalent in $\ZFC$ but not necessarily in the weaker framework of $\ZF$. This is particularly important when working in the context of the Axiom of Determinacy as the Axiom of Choice fails in this context. When (in contrast to Definition \ref{def:CH}) asking whether $|\bR| = |\omega_1|$, one could argue that the Continuum Hypothesis fails under the Axiom of Determinacy because $|\bR|$ and $|\omega_1|$ are incomparable. But this is not the path we want to follow in this article. Instead, we focus on Cantor's original formulation of the Continuum Problem, which is formalized in Definition \ref{def:CH}. See \cite{stanfordCH} for details on different formulations of $\CH$.
\end{remark}

Large cardinal axioms do not influence the Continuum Problem, i.e., it is still independent over $\ZFC$ together with the existence of large cardinals. The first instance of this phenomenon is given by the following result of Levy and Solovay on forcing in the context of measurable cardinals:

\begin{fact}[Levy-Solovay]
    If there is a measurable cardinal, there is a model with a measurable cardinal in which $\CH$ fails as well as a model with a measurable cardinal in which $\CH$ holds.
\end{fact}

The same result holds for all other known large cardinal axioms such as Woodin cardinals or supercompact cardinals. At the same time, sufficiently strong forcing axioms imply that the Continuum Hypothesis fails:

\begin{fact}
    Under the Proper Forcing Axioms there is a set $A$ such that  $|\mathbb{N}| < |A| < |\mathbb{R}|$. Even more, there is exactly one such intermediate size between  $|\mathbb{N}|$ and $|\mathbb{R}|$.
\end{fact}

This might seem like it is the end of the dream to solve independence phenomena in mathematics by considering canonical hierarchies of axioms extending $\ZF$. In fact, this is only the beginning of a long and fascinating line of research that we will partially outline in the rest of this article. 

\section{Identifying definable parts of the universe}

One important aspect of G\"{o}del's Program is the identification of the \textit{core} of the universe. Loosely speaking, the core of the universe consists of definable objects whose definability persists through all generic extensions of the universe. Examples of objects in the core are all the sets in G\"{o}del's constructible universe $L$ (see also Section \ref{sec:LandUltimateL}). The objects in the core that are of central interest are grounded in the universally Baire sets. The nature of such objects becomes visible in the context of the Axiom of Determinacy or Woodin's extension of $\AD$ called $\AD^+$. $\AD^+$ extends the Axiom of Determinacy by, in the context of the axiom of Dependent Choice for sets of reals ($\DC_\bR$), additionally asking that all sets of reals are $\infty$-Borel and ${<}\Theta$-ordinal determinacy holds. Whether $\AD$ implies $\AD^+$ is a central open problem in the area, it does so in all known natural models of determinacy. In this article, the reader unfamiliar with $\AD^+$ theory will lose little by ignoring the $+$. For a thorough exposition of $\AD^+$ theory we recommend \cite{La23}.

Cantor himself already showed that closed sets of reals cannot provide counterexamples to the Continuum Hypothesis:

\begin{theorem}[Cantor]
    Every closed set of reals is either countable or contains a perfect subset. Therefore, the Continuum Hypothesis holds for closed sets.
\end{theorem}

Indeed, a strengthening of Cantor's result holds: Every analytic set is either countable or contains a perfect subset. Therefore, the Continuum Hypothesis holds for analytic sets.

Large cardinals imply that the Continuum Hypothesis holds for larger classes of definable sets of reals. In fact, just like in Cantor's result, they imply that the classes of definable sets of reals in question have a stronger structural property, the perfect set property. This means, that the sets of reals in question are either countable or contain a perfect set. As perfect sets by nature have size continuum, the perfect set property implies the Continuum Hypothesis (in the version stated in Definition \ref{def:CH}). This implication holds in this general form but also when restricted to sufficiently closed classes of definable sets of reals. The reason why large cardinals imply that there are no definable counterexamples to $\CH$ in this setting is that they imply determinacy for the collections of sets of reals in question, which in turn implies the perfect set property. 

The first example (in terms of the complexity of the involved sets of reals) of a large cardinal implying a restricted version of the Continuum Hypothesis is the following:

\begin{theorem}[Solovay, \cite{So69}]
    Suppose there is a measurable cardinal. Then the Continuum Hypothesis holds for projections of analytic sets, i.e., for $\SIGMA^1_2$ sets.
\end{theorem}

In fact, Solovay showed directly that measurable cardinals imply the perfect set property for projections of analytic sets without going through determinacy. Alternatively, one can argue that, by a result of Martin \cite{Ma70}, measurable cardinals imply $\PI^1_1$ determinacy, which, in turn, implies the perfect set property for $\SIGMA^1_2$ sets.

In the context of stronger large cardinals, the Continuum Hypothesis holds for a larger collection of definable sets of reals.

\begin{theorem}[Steel, Woodin, \cite{St02, La04, St10}]
    Suppose there are infinitely many Woodin cardinals with a measurable cardinal above their supremum. Then the Continuum Hypothesis holds for all sets of reals constructible from $\bR$, i.e., all sets of reals in $L(\bR)$.
\end{theorem}

In fact, in a situation similar to the Steel-Woodin theorem a stronger result holds. It is not only the case that there are no constructible counterexamples to $\CH$ but, in addition, no such counterexamples can be added by (set-sized) forcing as the theory of $L(\bR)$ is \emph{frozen} (or \emph{sealed}) in the following sense:

\begin{theorem}[Steel, Woodin, \cite{St02, La04, St10}]\label{thm:freezingLR}
    Suppose there is a proper class of Woodin cardinals. Then the Axiom of Determinacy holds in $L(\bR)$ and for any pair of set generic extensions $V[g] \subseteq V[g*h]$ there is an elementary embedding \[ j \colon L(\bR_g) \rightarrow L(\bR_{g*h}), \] where $\bR_g$ and $\bR_{g*h}$ denote the sets of reals in $V[g]$ and $V[g*h]$ respectively.
\end{theorem}

A natural question at this point is whether it is possible to extend this result to sets of reals beyond $L(\bR)$. The first step in this direction is a result of Woodin extending Theorem \ref{thm:freezingLR} to $L(A,\bR)$ for an arbitrary universally Baire set $A$. Universally Baire sets originate in work of Schilling and Vaught \cite{SchVa83}, and they were first systematically studied by Feng, Magidor, and Woodin \cite{FMW92}. Since then they have played a prominent role in many areas of set theory. A set of reals is universally Baire if all of its continuous preimages in compact Hausdorff spaces have the property of Baire.  However, the following equivalent definition is set theoretically more useful. 

\begin{definition}[Feng-Magidor-Woodin, \cite{FMW92}]\label{def:uB}
\vspace{0.5em}\begin{enumerate}\itemsep0.5em
    \item  Let $(S,T)$ be trees on $\omega \times \kappa$ for some ordinal $\kappa$ and let $Z$ be any set. We say \emph{$(S,T)$ is $Z$-absolutely complementing} iff \[ p[S] = \BS \setminus p[T] \] in every $\Col(\omega,Z)$-generic extension of $V$.
    \item  A set of reals $A$ is \emph{universally Baire (uB)} if for every $Z$, there are $Z$-absolutely complementing trees $(S,T)$ with $p[S] = A$.
\end{enumerate}\vspace{0.5em}
\end{definition}

We write $\uB$ for the set of universally Baire sets and, if $g$ is set generic over $V$, we write $\uB_g = (\uB)^{V[g]}$ and $\bR_g = \bR^{V[g]}$ for the set of universally Baire sets and the set of reals in $V[g]$. Given a universally Baire set $A$ and a set generic $g\subseteq \mathbb{P}$ for some poset $\bP$, we write $A_g$ for the interpretation of $A$ in $V[g]$. More precisely, letting $\kappa\geq \card{\mathbb{P}}^+$ be a cardinal and $(S, T)$ be any $\kappa$-absolutely complementing trees with $p[S]=A$, we set $A_g=(p[S])^{V[g]}$. It can be easily checked using the absoluteness of well-foundedness that $A_g$ is independent of $(T, S)$. 

Assuming a proper class of Woodin cardinals, the study of the model $L(A, \mathbb{R})$, where $A$ is a universally Baire set, is completely parallel to the study of $L(\mathbb{R})$, and most of the theorems proven for $L(\mathbb{R})$ can be easily generalized to $L(A, \mathbb{R})$. For example, assuming a proper class of Woodin cardinals, generalizing the proof for $L(\mathbb{R})$, Woodin showed that the theory of the model $L(A, \mathbb{R})$ cannot be changed by forcing.

\begin{theorem}[Woodin]\label{thm:freezingLAR}
    Suppose there is a proper class of Woodin cardinals and let $A$ be an arbitrary universally Baire set of reals. Then the Axiom of Determinacy (even $\AD^+$) holds in $L(A,\bR)$ and for any pair of set generic extensions $V[g] \subseteq V[g*h]$ there is an elementary embedding \[ j \colon L(A_g, \bR_g) \rightarrow L(A_{g*h}, \bR_{g*h}) \] 
    such that $j\restriction \mathbb{R}_g=id$ and $j(A_g)=A_{g*h}$.
\end{theorem}

It is known that sets of reals that are definable via sufficiently generically absolute formulas are universally Baire (e.g., see \cite[Lemma 4.1]{St09}). In this sense, the generically absolute fragment of $\ZFC$ is coded into the universally Baire sets, and in this sense, the model $L(\uB, \bR)$, as far as the sets of reals go, must be the maximal generically absolute inner model. However, the story of $L(\uB, \bR)$ is not as simple as one might guess.
Extending Theorem \ref{thm:freezingLAR} from a single universally Baire set $A$ to the set of all universally Baire sets of reals $\uB$ turns out to be non-trivial and is, in fact, a deep question at the surface of our current knowledge in set theory. We will make the statement more precise and summarize the state of the art on this question in the next section.

\section{Sealing the universally Baire sets and beyond}\label{sec:Sealing}

In contrast to the situation for $L(A, \bR)$ for a single universally Baire set $A$, it is not clear that assuming large cardinals, $L(\uB, \bR)\models \AD^+$ or even 
\begin{center}
    $\powerset(\mathbb{R})\cap L(\uB, \bR)=\uB$.
\end{center}

Woodin isolated the following notion which states the natural extension of Theorem \ref{thm:freezingLAR} to the set of all universally Baire sets of reals. Informally, it states that the theory of the universally Baire sets is \emph{sealed}.

\begin{definition}[Woodin]\label{def:Sealing}
$\Sealing$ is the conjunction of the following statements.
\vspace{0.5em}\begin{enumerate}\itemsep0.5em
    \item For every set generic $g$ over $V$, $L(\uB_g, \bR_g) \models \AD^+$ and $\powerset(\bR_g) \cap L(\uB_g, \bR_g) = \uB_g$. \label{eq:Sealing1}
    \item For every set generic $g$ over $V$ and set generic $h$ over $V[g]$, there is an elementary embedding \[ j \colon L(\uB_g, \bR_g) \rightarrow L(\uB_{g*h}, \bR_{g*h}) \] such that for every $A \in \uB_g$, $j(A) = A_h$. \label{eq:Sealing2}
\end{enumerate}\vspace{0.5em}
\end{definition}

Using the stationary tower forcing, Woodin showed that $\Sealing$ is consistent from a supercompact cardinal and a proper class of Woodin cardinals. 

\begin{theorem}[Woodin's Sealing Theorem, see \cite{La04,WoodinLongExtender}]\label{thm:SealingWoodin}
Let $\kappa$ be a supercompact cardinal and let $g$ be $\Col(\omega, 2^\kappa)$-generic over $V$. Suppose that there is a proper class of Woodin cardinals. Then $\Sealing$ holds in $V[g]$.
\end{theorem}

Woodin, in addition, showed that in this setting $\Theta$ is regular in $L(\uB_g, \bR_g)$, see \cite{WoodinLongExtender, MuSa_ThetaReg}. Here $\Theta$ denotes the supremum of all ordinals onto which the reals $\bR$ surject. Theorem \ref{thm:SealingWoodin} implies that in $V[g]$, $L(\uB, \bR)$ is a canonical model whose theory cannot be changed by forcing. The authors have shown in \cite{MuSa_Sealing} that the reason for this is that $L(\uB, \bR)$ can be realized as a derived model at a supercompact cardinal via Woodin's Derived Model Theorem (see Theorem \ref{woodin: der model thm}). Recall that for any $Z,X$ and any ordinal $\gamma$, $\meas_\gamma(Z)$ denotes the set of all $\gamma$-additive measures on $Z^{{<}\omega}$. 
We write $\bar\mu = (\mu_s \mid s \in X^{{<}\omega})$ for a $\gamma$-complete \emph{homogeneity system} over $X$ with support $Z$ 
if for each $s \in X^{{<}\omega}$, $\mu_s \in \meas_\gamma(Z)$. For details on the definition of homogeneity systems we refer the reader to \cite{St09}.

\begin{definition}\label{def: smu}
    A set $A \subseteq X^{\omega}$ is $\gamma$-homogeneously Suslin if there is a $\gamma$-complete homogeneity system $\bar\mu = (\mu_s \mid s \in X^{{<}\omega})$ and a tree $T$ such that $A = p[T]$ and, for all $s \in X^{{<}\omega}$, $\mu_s(T_s) = 1$. In particular, \[ A = S_{\bar{\mu}}=_{def} \{ x \in X^\omega \mid (\mu_{x \upharpoonright n} \mid n < \omega) \text{ is well-founded} \}.\footnote{This means that the natural direct limit of $\Ult(V, \mu_{x \upharpoonright n})$ for $n<\omega$ is well-founded.} \]
\end{definition}

Here $T_s = \{ t \mid (s, t) \in T \}$. We write \[ \Hom_\gamma = \{ A \subseteq X^{\omega} \mid A \text{ is } \gamma\text{-homogeneously Suslin} \} \] and \[ \Hom_{{<}\eta} = \bigcup_{\gamma < \eta} \Hom_\gamma. \]

For $\eta$ a limit of Woodin cardinals and $H \subseteq \Col(\omega, {<}\eta)$ generic over $V$, write $H \upharpoonright \alpha = H \cap \Col(\omega, {<}\alpha)$. As usual, let \[ \bR^* = \bR_H^* = \bigcup_{\alpha < \eta} \bR \cap V[H \upharpoonright \alpha]. \] Moreover, for any $\alpha < \eta$ and $A \in \Hom_{{<}\eta}^{V[H \upharpoonright \alpha]}$, write \[ A^* = \bigcup_{\alpha < \beta < \eta} A_{H \upharpoonright \beta} \] and \[ \Hom^* = \{ A^* \mid \exists \alpha < \eta \, A \in (\Hom_{{<}\eta})^{V[H\upharpoonright \alpha]}\}. \]

In this setting, $L(\bR^*, \Hom^*)$ is called a derived model of $V$ at $\eta$. The following is Woodin's Derived Model Theorem (see \cite[Theorem 6.1]{St09}).

\begin{theorem}[Woodin, see \cite{St09, StstattowfreeDMT}]\label{woodin: der model thm} Suppose $\eta$ is a limit of Woodin cardinals and $H \subseteq \Col(\omega, {<}\eta)$ generic over $V$. Let $L(\bR^*, \Hom^*)$ be the model computed as above in $V(\bR^*)$. Then $L(\bR^*, \Hom^*) \models \AD^+$.
\end{theorem}

The authors proved the following theorem which implies a version of Theorem \ref{thm:SealingWoodin} when the generic $g$ is chosen as $\Col(\omega, 2^{2^\kappa})$-generic (instead of just $\Col(\omega, 2^\kappa)$-generic). The argument was inspired by earlier work of Trang and the second author on $\Sealing$ from iterability \cite{SaTr21}.

\begin{theorem}[Müller-Sargsyan, \cite{MuSa_Sealing}]\label{thm:MuellerSargsyanSealing}
    Let $\kappa$ be a supercompact cardinal and suppose there is a proper class of Woodin cardinals. Let $g \subseteq \Col(\omega, 2^\kappa)$ be $V$-generic, $h$ be $V[g]$-generic and $k\subseteq \Col(\omega, 2^\omega)$ be $V[g*h]$-generic. Then, in $V[g*h*k]$, there is $j: V\rightarrow M$ such that $j(\kappa)=\omega_1^{M[g*h]}$ and  $L(\uB_{g*h}, \bR_{g*h})$ is a derived model of $M$, i.e., for some $M$-generic $G \subseteq \Col(\omega, {<}\omega_1^{V[g*h]})$,
        \[ L(\uB_{g*h}, \bR_{g*h}) = (L(\Hom^*, \bR^*))^{M[G]}. \]
\end{theorem}

A key question is whether these results can be extended to larger parts of the universe:

\begin{question}
 What sort of canonical, generically absolute subsets of $\uB$ can be added to the model $L(\uB, \bR)$ in the context of $\Sealing$?
\end{question}

The authors consider extensions of $L(\uB, \bR)$ in two different directions in \cite{MuSa_Sealing}: by adding the club filter and by adding the ordinal definable powerset of $\uB$. We focus on the definable powerset of $\uB$ here and refer the interested reader to \cite{MuSa_Sealing} for the results on $\Sealing$ with the club filter on $\powerset_{\omega_1}(\uB)$.

The definable powerset of $\uB$ that we add to $L(\uB, \bR)$ is defined as follows.

\begin{definition}\label{ub powerset} Suppose $X$ is a set. We let $\iota_X=\max(\card{X}, \card{\uB})$ and $\powerset_{uB}(X)$, the \emph{$\uB$-powerset of $X$}, be the set of those $Y\subseteq X$ such that whenever $g\subseteq \Col(\omega, \iota_X)$ is $V$-generic, $Y$ is ordinal definable in $L(\uB_g, \bR_g)$ from parameters belonging to the set $\{X, j_g \shortpwimg \uB \}\cup j_g \shortpwimg \uB$, where $j_g: L(\uB, \bR)\rightarrow L(\uB_g, \bR_g)$ is the canonical embedding with the property that $j_g(A)=A_g$ for every  $A \in \uB$.   
\end{definition}

 The following is an easy consequence of $\sf{Sealing}$. The proof is outlined in the introduction of \cite{MuSa_Sealing}.
    \begin{lemma}\label{easy consequence} Assume $\sf{Sealing}$ and suppose $Y\in \powerset_{uB}(X)$. Suppose $h$ is any $V$-generic with the property that there is $k\in V[h]$ which is $V$-generic for $\Col(\omega, \iota_X)$. Then $Y$ is ordinal definable in $L(\uB_h, \bR_h)$ from parameters that belong to the set $\{X, j_h \shortpwimg \uB \}\cap j_h \shortpwimg \uB$.\footnote{Here $j_h \shortpwimg \uB$ denotes the pointwise image of $\uB$ under $j_h$.}
    \end{lemma}
   
Let $\uBp=\powerset_{uB}(\uB)$. If $g$ is $V$-generic then we let $\uBp_g=(\uBp)^{V[g]}$. If $\sf{Sealing}$ holds and $g$, $g'$ are two consecutive generics\footnote{I.e. $g'$ is $V[g]$-generic.}, then we let
\begin{center}
    $j_{g, g'}:L(\uB_g, \bR_g)\rightarrow L(\uB_{g*g'}, \bR_{g*g'})$
\end{center}
be the canonical $\sf{Sealing}$ embedding, i.e., $j_{g, g'}(A)=A_{g'}$.

\begin{definition}\label{def: sealing for the ub powerset} We say $\sf{Weak\ Sealing}$ holds for the uB-powerset\footnote{This is different from the Weak Sealing notion in \cite{WoodinLongExtender}.} if 
\vspace{0.5em}\begin{enumerate}\itemsep0.5em
    \item $\Sealing$ holds, \label{eq: sealing for the ub powerset clause 1}
    \item $L(\uBp)\models \AD^+$, \label{eq: sealing for the ub powerset clause 2}
    \item whenever $g, g'$ are two consecutive generics such that $V[g*g']\models \card{(2^{2^\omega})^{V[g]}}=\aleph_0$, there is an elementary embedding $\pi:L(\uBp_g)\rightarrow L(\uBp_{g*g'})$ such that $\pi\restriction L(\uB_g, \bR_g)=j_{g, g'}$. \label{eq: sealing for the ub powerset clause 3}
\end{enumerate}\vspace{0.5em}
We say $\sf{Sealing}$ holds for the uB-powerset if Clause \eqref{eq: sealing for the ub powerset clause 3} above holds for all consecutive generics $g$ and $g'$.
\end{definition}
\begin{remark}
If $\sf{Weak\ Sealing}$ holds for the uB-powerset then the theory of the model $L(\uBp)$ cannot be changed by forcing. This is because if $g$ and $g'$ are two consecutive generics and $W$ is a generic extension of both $V[g]$ and $V[g*g']$ such that $W\models \card{\powerset(\powerset(\bR))^{V[g]}}=\card{\powerset(\powerset(\bR))^{V[g*g']}}=\aleph_0$ then both $L(\uBp_g)$ and $L(\uBp_{g*g'})$ are elementarily equivalent to $L(\uBp)^W$.
\end{remark}

\begin{theorem}[Müller-Sargsyan, \cite{MuSa_Sealing}]\label{sealing for the ub powerset} Suppose $\kappa$ is a supercompact cardinal and there is a proper class of inaccessible limits of Woodin cardinals. Suppose $g\subseteq \Col(\omega, 2^{2^\kappa})$ is $V$-generic. Then $\sf{Weak\ Sealing}$ for the uB-powerset holds in $V[g]$. 
\end{theorem}

Theorem \ref{sealing for the ub powerset} implies that $L(\uBp)$ is a canonical model of determinacy. Moreover, in the setting of Theorem \ref{sealing for the ub powerset}, $\Theta$ is regular in $L(\uBp)$ by the results in \cite{MuSa_ThetaReg}.

\begin{question} Does ${\sf{Sealing}}$ for the uB-powerset hold in the model of Theorem \ref{sealing for the ub powerset} or in any other model? 
\end{question}

\begin{remark} We conjecture that $\uBp$ is the amenable powerset of $\uB$. More precisely, $\uBp$ is the set of all $A\subseteq \uB$ such that for every $\k<\Theta^{L(\uB, \bR)}$, letting $\Gamma_\k=\{ B\in \uB \mid$ the Wadge rank of $B$ is $<\k\}$, $A\cap \Gamma_\k\in L(\uB, \bR)$. This conjecture is closely related to Conjecture \ref{cofinality conjecture} below.
\end{remark}

A key motivation for considering these extensions of $L(\uB, \bR)$ in the context of $\Sealing$ is the possibility to force ${\sf{MM^{++}}}$ over such determinacy models as well as their siblings. This has been a major line of set theoretic research in the past years and we refer the interested reader to the introduction of \cite{MuSa_Sealing} as well as to \cite{SaAnnouncement, LaSa21} for details. To follow this line of research, we need to extend the model $L(\uBp)$ even further as it is unlikely that ${\sf{MM^{++}}(c^{++})}$ can be forced over $L(\uBp)$. The further extension is motivated by the model used by Larson and the second author in \cite{LaSa21}. The next step is to consider the model $L(\powerset_{\omega_1}(\uBp))$ and prove a ${\sf{Sealing}}$ theorem for it. For the model to make sense, however, we need to know that $\sup (Ord\cap \uBp)$ has cofinality $\omega$. 
 
\begin{definition}\label{thetax} Suppose $X$ is a set. Let $\eta_X$ be the supremum of all ordinals $\a$ such that there is a surjective $f \colon \powerset(\bR) \rightarrow \a$ which is ordinal definable from $X$.
\end{definition}

\begin{definition}\label{def: wst} We say that the $\sf{Weak\ Sealing\ Theorem}$ for a formula $\phi$ holds if setting $M_\phi=\{ x: V\models \phi[x]\}$, the following statements hold:
\vspace{0.5em}\begin{enumerate}\itemsep0.5em
\item $M_\phi$ is a transitive model of ${\sf{AD^+}}$.
\item For all $V$-generics $g$, the models $(M_\phi)^{V}$ and $(M_\phi)^{V[g]}$ are elementarily equivalent.
\end{enumerate}\vspace{0.5em} 
We also simply say that the $\sf{Weak}$ ${\sf{Sealing}}$ ${\sf{Theorem}}$ for $M_\phi$ holds to mean that the $\sf{Weak\ Sealing\ Theorem}$ holds for $\phi$. Similarly, we define the meaning of ${\sf{Sealing}}$ $\sf{Theorem}$ for $\phi$. This means that  
\vspace{0.5em}\begin{enumerate}\itemsep0.5em
\item $M_\phi$ is a transitive model of ${\sf{AD^+}}$ containing all the reals.
\item for all $V$-generic $g$, $V[g]\models \powerset(\bR)\cap M_\phi \subseteq \uB$.
\item For all $V$-generic $g$ and $V[g]$-generic $h$, there is an elementary embedding $j:(M_\phi)^{V[g]}\rightarrow (M_\phi)^{V[g*h]}$ such that for every $A\in \powerset(\bR)^{(M_\phi)^{V[g]}}$, $j(A)=A_h$.\label{it:Mphi3}
\end{enumerate}\vspace{0.5em} 
We say that the Sealing embedding for $\phi$ preserves $\psi$ if, letting $B_\psi=\{a: V\models \psi[a]\}$, 
\begin{enumerate}\itemsep0.5em
\item for every $V$-generic $g$, $(B_\psi)^{V[g]}\cap (M_\phi)^{V[g]} \in (M_\phi)^{V[g]}$ and
\item there exists $j$ as in Item \eqref{it:Mphi3} above with the additional property that $$j((B_\psi)^{V[g]}\cap (M_\phi)^{V[g]})=(B_\psi)^{V[g*h]}\cap (M_\phi)^{V[g*h]}.$$
\end{enumerate}\vspace{0.5em} 
\end{definition}

\begin{definition}\label{changmartin} Let $\eta = (\eta_{\uB})^{L(\uB, \bR)}$ and let $\mathsf{C-uBp}=L(\eta^\omega, \uBp)$ be the \emph{Chang-type $\uB$-powerset model}. 
\end{definition}

\begin{conjecture}\label{cofinality conjecture} Suppose $\kappa$ is a supercompact cardinal and there is a proper class of inaccessible limits of Woodin cardinals. Suppose $g\subseteq \Col(\omega, 2^{2^\kappa})$ is $V$-generic. Let $\eta_g= (\eta_{\uB_g})^{L(\uB_g, \bR_g)}$. Then $\cf(\eta_g)=\omega$. Moreover, the ${\sf{Sealing}}$ ${\sf{Theorem}}$ holds for $\mathsf{C-uBp}$ in $V[g]$.
\end{conjecture}

\section{Canonical inner models and Woodin's axiom $\UltimateL$}\label{sec:LandUltimateL}

As studying the mathematical universe from above by isolating its global properties can be difficult (if not even impossible), a natural approach is to study the universe bottom-up, starting from the parts that are easier to analyze. This is the approach Gödel took when he defined the Constructible Universe $L$. In fact, it is more accurate to call this approach \emph{inside-out} as $L$ is as high as the full universe $V$ but it will (for example, in the context of measurable cardinals) be thinner than $V$. For the readers' convenience we briefly recall the definition of Gödel's Constructible Universe $L$ and its basic properties. $L$ is constructed bottom-up inductively throughout the ordinals and the successor steps will be given by the following operation:

\begin{definition}
  For any set $A$, let $\Def(A)$ be the set of all subsets of $A$ that are definable over $(A,\in)$ with parameters in $A$. Here we say that $X \subseteq A$ is \emph{definable over $(A,\in)$} iff there is a formula $\varphi$ (in the language of set theory with one relational symbol $\in$) and some $a_1, \dots, a_n \in A$ such that \[ X = \{ x \in A \mid (A, \in) \models \varphi[x,a_1,\dots,a_n]\}. \]
\end{definition}

Note that $A \in \Def(A)$ and, if $A$ is transitive, $A \subseteq \Def(A) \subseteq \powerset(A)$. Now we are ready to define G\"odel's Constructible Universe $L$.

\begin{definition}[Gödel]
  We define via transfinite induction:
  \begin{eqnarray*}
    L_0 &=& \emptyset,\\
    L_{\alpha+1} &=& \Def(L_\alpha),\\
    L_\lambda &=& \bigcup_{\alpha<\lambda} L_\alpha, \text{ if } \lambda \text{ is a limit ordinal, and}\\
    L &=& \bigcup_{\alpha \in \Ord} L_\alpha.
  \end{eqnarray*}
  The proper class $L$ is \emph{Gödel's Constructible Universe}.
\end{definition}

For each level $\alpha$, $L_\alpha \subseteq V_\alpha$ and the ordinal height of $L_\alpha$ is $\alpha$, i.e., $L_\alpha \cap \Ord = \alpha$, and in particular $\alpha \subseteq L_\alpha$. Moreover, $|L_\alpha| = |\alpha|$ for all $\alpha \geq \omega$. The key properties of $L$ are that the Axiom of Choice $\AC$ and the Continuum Hypothesis $\CH$ (even the Generalized Continuum Hypothesis $\GCH$) hold in $L$. In fact, we can even show this more generally using the Axiom of Constructibility. The \emph{Axiom of Constructibility} ``$V=L$'' states that $\forall x \exists \alpha (x \in L_\alpha)$. By absoluteness, ``$V = L$'' holds in $L$. Moreover, again using absoluteness, it is not difficult to show that, given $\ZF$ in the universe $V$, $\ZF$ holds in $L$. Now the following theorem yields $\AC$ and $\GCH$ in the Constructible Universe $L$:

\begin{theorem}[Gödel]
  $\ZF$ + ``$V = L$'' proves $\AC$ and $\GCH$.
\end{theorem}

The fact that ``$V = L$'' proves the Axiom of Choice follows from the observation that every $L_\alpha$ can be well-ordered. First, $\Def(A)$ can be canonically well-ordered given a well-order of $A$. Then we piece these well-orders together to get a well-order of $L_\alpha$. This, in fact, shows more than just the Axiom of Choice. It shows that ``$V=L$'' implies \emph{global choice}, in fact, there is a definable relation $<_L$ that well-orders all sets in $L$.

The Generalized Continuum Hypothesis $\GCH$ in $L$ is a consequence of the \emph{Condensation Lemma}, a key property of G\"odel's constructible universe $L$. It can be used to show that ``$V=L$'' implies that $L_\kappa = H(\kappa)$ for all cardinals $\kappa \geq \omega$, where $H(\kappa) = \{ x \in \WF \mid |\trcl(x)| < \kappa\}$ for $\WF$ the class of all well-founded sets and $\trcl(x)$ is the transitive closure of $x$. This together with ``$V=L$'' implies $2^\kappa = \kappa^+$ for all cardinals $\kappa \geq \omega$ and hence $\GCH$.

These results demonstrate that from a certain perspective Gödel's Constructible Universe $L$ is a very desirable model: it is well understood and satisfies properties such as the (Generalized) Continuum Hypothesis. Unfortunately, it is less useful when combined with the large cardinal hierarchy: already mild large cardinals cannot exist in $L$.

\begin{theorem}[Scott \cite{Sc61}]
    $\ZF$ + ``$V = L$'' proves that there are no measurable cardinals.
\end{theorem}

But this is not the end of the story of Gödel's Constructible Universe, it is in fact the birth of an active and flourishing research area called inner model theory aiming to follow Gödel's approach while incorporating stronger and stronger large cardinals into the models. To illustrate the main ideas as well as obstacles, we sketch the construction of an $L$-like model with a measurable cardinal here.

Recall that a \emph{measure} $U$ on an uncountable cardinal $\kappa$ is a $\kappa$-complete non-principal ultrafilter $U$ on $\kappa$. It is possible to construct $L$ relative to the predicate $U$\footnote{That means allowing $U$ as a predicate in the language for $\Def$ in the successor steps of the construction.} and obtain a $\ZFC$ model $L[U]$ with a measurable cardinal.

\begin{theorem}[Solovay]
    Let $U$ be a measure on some cardinal $\kappa$. Then $U \cap L[U] \in L[U]$ and $U \cap L[U]$ is a measure on $\kappa$. 
\end{theorem}

The Condensation Lemma for $L[U]$ is more involved but the following result shows that $L[U]$ indeed resembles Gödel's Constructible Universe $L$.

\begin{theorem}[Kunen, Silver] 
Let $U$ be a measure on some cardinal $\kappa$ in $L[U]$.
\begin{enumerate}
    \item Let $U'$ be a measure on the same cardinal $\kappa$ in $L[U']$. Then $L[U] = L[U']$ and $U \cap L[U] = U' \cap L[U]$.
    \item \GCH holds in $L[U]$.
\end{enumerate} 
\end{theorem}

A key method used in proving the above theorem is the theory of iterated ultrapowers going back to Gaifman. Two models of the same form, say $L[U]$ and $L[U']$ can be compared by considering linear iterations of consecutively applying their ultrafilters and their images, taking direct limits at limit stages of this process.

Solovay showed that models of the form $L[U]$ have a unique measurable cardinal. So altogether they are \emph{canonical inner models for a measurable cardinal} but, analogously to Scott's Theorem for Gödel's $L$, they cannot accommodate any stronger large cardinals than measurables. This started a long line of research to construct inner models for stronger large cardinals. For many measurable cardinals and strong cardinals this work was pioneered by Ronald B. Jensen, Robert M. Solovay, Tony Dodd and William J. Mitchell and it was extended to the level of Woodin cardinals by Donald A. Martin, William J. Mitchell and John R. Steel \cite{MaSt94, MS94, SchStZe02}. To build canonical inner models with strong cardinals or Woodin cardinals, instead of constructing relative to a single measure or extender,\footnote{An extender is a system of measures that allows to build stronger ultrapowers which can capture arbitrary elementary embeddings.} the model $L[\vec E]$ is built relative to a coherent sequence of extenders $\vec E$. Proving that models of the form $L[\vec E]$ resemble $L$ and, for example satisfy the \GCH, requires a very fine analysis of their internal structure. This was pioneered by Jensen who developed a \emph{fine structure analysis}  and also introduced the terminology \emph{mouse} for such canonical models \cite{Je72}. In fact, what we just described is formally known as a \emph{premouse}. See, for example, \cite{St10} for the full definition.

\begin{definition}\label{def:premouse}
    A structure of the form $(L_\alpha[\vec E], \vec E, F)$ for $\alpha \leq \OR$ is called a \emph{premouse} if
    \begin{enumerate}
        \item $\vec E$ is a coherent sequence of extenders with a canonical indexing of the extenders, 
        \item $F$ is an extender over $L_\alpha[\vec E]$ with index $\alpha$ or $F = \emptyset$,
        \item $\vec E$ is \emph{acceptable} at each $\beta < \alpha$, meaning that for every $\kappa$, if $L_{\beta+1}[\vec E]$ has a subset of $\kappa$ that is not an element of $L_\beta[\vec E]$, then $|L_\beta[\vec E]| \leq \kappa$ holds in $L_{\beta+1}[\vec E]$, and
        \item all proper initial segments of $(L_\alpha[\vec E], \vec E, F)$ satisfy a certain fine structural definability condition called \emph{soundness}.\index{Soundness}
    \end{enumerate}
\end{definition}

Note that the condition that $\vec E$ is acceptable at each $\beta < \alpha$ implies that the Continuum Hypothesis, and more generally $\GCH$, is true in premice. To analyze these models we again need a method of comparison by taking iterated ultrapowers. This turns out to be more complicated than before for inner models with Woodin cardinals. More precisely, these iterations need to be organized as trees, allowing us to use extenders from the current model on earlier models on different branches of the tree as long as they agree sufficiently far for this to make sense. This leaves us with the question of how to proceed at limit stages. In the setting above for $L[U]$ the iteration was linear and at a limit stage the only option was to take the direct limit along this linear iteration. But an infinite tree can have more than one cofinal branch and it becomes a non-trivial question which well-founded model arousing as a direct limit along a cofinal branch through our tree we should choose as the next model of our iteration tree. Therefore, the canonical inner models need iteration strategies which tell us along which branch of an iteration tree we take the direct limit at a limit stage of the iteration.
It turns out that to build canonical inner models for large cardinals at the level of Woodin cardinals and beyond, iteration strategies and their complexity are the key property of mice.

This method of building canonical inner models for large cardinals is very successful for a certain part of the large cardinal hierarchy, currently roughly up to the level of a Woodin limit of Woodin cardinals. When considering stronger large cardinal notions or aiming to extend this all the way, the natural question emerges whether the set theoretic universe $V$ itself is a canonical model. If so, it would be the ultimate version of Gödel's constructible universe $L$. Woodin isolated the following formalization of this setting. 

\begin{definition}[Woodin, $\UltimateL$]\hspace{0.5cm}\\
The axiom $\UltimateL$ is the conjunction of the following statements:
      \begin{enumerate}
      \item There is a proper class of Woodin cardinals.
      \item For each $\Sigma_2$-sentence $\phi$, if $\phi$ holds in $V$ then there exists a universally Baire set $A \subseteq \mathbb{R}$ such that \[ \HOD^{L(A,\mathbb{R})} \models \phi. \]
      \end{enumerate}
\end{definition}

The motivation behind $\UltimateL$ is based on an extensive study of $\HOD$ in models of determinacy called $\HOD$ analysis. This study suggests that in the context of determinacy $\HOD$, the class of all hereditarily ordinal definable sets, is indeed a canonical inner model but not quite a premouse as introduced in Definition \ref{def:premouse}. It is a \emph{strategy mouse}, that means, in addition to an extender sequence it has a predicate for its own iteration strategy. Before we describe this in a bit more detail, note that the Continuum Hypothesis holds if $\UltimateL$ for the following reason. Both $\CH$ and $\neg\CH$ can be written as $\Sigma_2$-sentences. As by $\AD^+$ theory it is known that $\CH$ holds in $\HOD^{L(A,\mathbb{R})}$ for universally Baire sets $A$ in the context of a proper class of Woodin cardinals, $\UltimateL$ naturally implies Cantor's Continuum Hypothesis. 

The analysis of $\HOD$ of models of determinacy is one of the most central areas of modern inner model theory. The program has its roots in the work of the Cabal group, in particular in the work of Becker and Moschovakis, but it firmly became a subject within inner model theory through the groundbreaking work of Steel and Woodin. Steel's seminal \cite{St95} initiated the study, and further investigations \cite{StW16} cemented the role of $\HOD$ of determinacy models in inner model theory and set theory in general. In addition to its role in Woodin's $\sf{Ultimate\ L}$ framework, $\HOD$ analysis played a key role in identifying mantles of canonical universes (see \cite{SaSch18}), it played a fundamental role in a recent solution of an old conjecture of Kechris (see \cite{sargsyan2021hjorths}), it plays a fundamental role in one of the most advanced inner model theoretic methods, the core model induction (see the discussion in Section \ref{sec:caninnermodelsSealing} and, for example, \cite{Mu20rev, Ketchersid, St05, Sa15covuB, TrWi21, ASTWZ24, SaTrConSealing}), it played a fundamental role in the refutation of the iterability conjecture for $K^c$  by Larson and the second author in \cite{LaSa21}, and finally it plays a fundamental role in establishing equiconsistencies between large cardinals and determinacy axioms, and in particular, in the recent solution of Sargsyan's Conjecture by the first author in \cite{Mu21}. 

The goal of $\HOD$ analysis is simply to show that $\HOD$ of models of determinacy is a canonical inner model. It naturally splits into two problems. 

\smallskip

\noindent\textbf{First $\HOD$ Analysis}:\\ Show that $V_\Theta^{\sf{HOD}}$ is a canonical fine structural inner model (and hence satisfies $\GCH$).  

\smallskip

\noindent\textbf{Second $\HOD$ Analysis}:\\ Show that $\HOD$ is a canonical fine structural inner model (and hence satisfies $\GCH$).

\smallskip

Steel's \cite{St95} completes the first $\HOD$ analysis assuming $V=L(\bR)+\AD$, and subsequent work of Woodin \cite{StW16} completes the second $\HOD$ analysis assuming $V=L(\bR)+\AD$. The division is natural because, as demonstrated by Woodin, $\HOD$ is fine structurally different than $V_\Theta^{\HOD}$. For example, in $L(\bR)$, $V_\Theta^{\HOD}$ is just an ordinary extender model of the form $L_\Theta[\vec{E}]$ while $\HOD$ is a hybrid extender model of the form $L[V_\Theta^{\HOD}, \Sigma]$ where $\Sigma$ is an iteration strategy for $V_\Theta^{\HOD}$. 

In \cite{Sa15}, Sargsyan completed the first $\HOD$ analysis while working in the minimal model of $\AD_\bR + \Theta \text{ is regular}$, and in \cite{Tr14}, Nam Trang completed the second $\HOD$ analysis while working in the same model or smaller ones. Currently, the best known result on the first and the second $\HOD$ analysis is the calculation of $\HOD$ of the minimal model of the \emph{Largest Suslin Axiom} in \cite{SaTr}. Steel's \cite{St22} develops general techniques for solving the first $\HOD$ analysis, and does so provided \emph{Hod Pair Capturing} holds. Thus, currently, proving Hod Pair Capturing is the central open problem in the area.


In the context of $\Sealing$ there is a canonical model of determinacy containing all universally Baire sets: $L(\uB, \bR)$. It was observed by Woodin that this can be used to define a canonical strengthening of $\UltimateL$. This turns $\UltimateL$ into an even more elegant and natural axiom and can arguably be seen as evidence for $\UltimateL$ and hence for $\CH$.

\begin{fact}[Woodin]
 In the context of $\Sealing$, the conjunction of the following two clauses is a well-defined strengthening of $\UltimateL$:
      \begin{enumerate}
      \item There is a proper class of Woodin cardinals.
      \item For each $\Sigma_2$-sentence $\phi$, if $\phi$ holds in $V$ then \[ \HOD^{L(\operatorname{uB},\mathbb{R})} \models \phi. \]
      \end{enumerate}
\end{fact}

\section{Canonical inner models and $\Sealing$}\label{sec:caninnermodelsSealing}

We have seen in Section \ref{sec:Sealing} that $\Sealing$ and several strengthenings of it can be obtained in generic extensions of the universe collapsing large cardinals in the region of supercompact cardinals. From an inner model theoretic perspective the following is a central question:

\begin{question}
    Is there a large cardinal that implies $\Sealing$?
\end{question}

Let us explain why. Let $M$ be a canonical inner model\footnote{For simplicity, let us assume $M$ is a proper class model. In fact, only very simple closure properties suffice.} with some large cardinals of the type described in Section \ref{sec:LandUltimateL}. Then, in $M$, there is a well-ordering of the reals in $L(\uB,\bR)$. But the first clause of $\Sealing$ for a trivial generic $g$ states that $L(\uB,\bR) \models \AD$. Since under $\AD$ there is no well-ordering of the reals, this yields that $\Sealing$ cannot hold in $M$. Therefore, the inner model theory of a large cardinal implying $\Sealing$ must be utterly different than the modern inner model theory as developed in \cite{St10} or other publications. This can be summarized as the well-known $\Sealing$ Dichotomy.

\begin{namedthm}{$\Sealing$ Dichotomy}
    Either no large cardinal notion implies $\Sealing$ or there is a large cardinal notion for which there is no canonical inner model of the type we currently know.
\end{namedthm}

While the known results on $\Sealing$ in generic extensions collapsing large cardinals outlined in Section \ref{sec:Sealing} have no direct consequence for constructing canonical inner models with large cardinals, there are results on $\Sealing$ that do have a direct impact on inner model theory. The following surprising result of Trang and the second author using deep inner model theoretic machinery shows that $\Sealing$ consistently appears well below a supercompact cardinal.

\begin{theorem}[Sargsyan-Trang, \cite{SaTrConSealing}]\label{thm:ConSealing}
    $\Sealing$ is consistent from a Woodin cardinal that is itself a limit of Woodin cardinals.
\end{theorem}

This result is no only surprising, it makes us change our current conception of consistency proofs via inner model theoretic methods. To explain this, we need to say a few words about the Core Model Induction, a technique to climb the consistency hierarchy via models of determinacy invented by Woodin in the 1990s and further developed by many inner model theorists since then.

To outline the core model induction technique we pick an example and summarize how consistency strength in terms of large cardinals can be obtained from the Proper Forcing Axiom $\PFA$ (or, technically, from the non-existence of certain square sequences). The starting point of the core model induction is Jensen's Covering Lemma for $L$ which yields, from $\PFA$, that $0^\sharp$ exists. By using a covering result for Steel's core model $K$, Schimmerling in \cite{Schim95} extended this result and showed that, under $\PFA$, $M_1^\sharp(X)$ exists for all sets $X$.\footnote{$M_1^\sharp(X)$ is the canonical sufficiently iterable inner model with one Woodin cardinal and a top extender build over $X$.} Steel and Woodin  independently extended this result to the existence of  $M_n^\sharp(X)$ for all natural numbers $n$ and all sets $X$, the canonical model for $n$ Woodin cardinals and a top extender. As a corollary, $\PFA$ implies projective determinacy. Steel's and Woodin's argument, inductively on $n$, uses a covering result for Steel's core model with finitely many Woodin cardinals. This already explains the terminology \emph{core model induction}. But the key ideas behind this methodology in fact only show up beyond the projective hierarchy when the focus shifts from inner models with large cardinals towards models of determinacy. The first step was Woodin's argument showing that $\PFA$ together with a strongly inaccessible cardinal yields the Axiom of Determinacy in $L(\bR)$. It requires a careful induction through the $L(\bR)$-hierarchy using descriptive set theoretic machinery to organize the induction. The first core model induction beyond $L(\bR)$ was done by Ketchersid in his PhD thesis \cite{Ketchersid}. 

The core model induction beyond $L(\bR)$ is the setting that is most relevant for us to draw the connection to $\Sealing$. In this setting the core model induction is arguably neither an induction nor a core model argument, it is a machinery for obtaining \emph{maximal models of determinacy} from various hypotheses.\footnote{These hypotheses will not be the existence of certain large cardinals but, for example, combinatorial statements or forcing axioms. The key method for obtaining models of determinacy from large cardinal axioms is Woodin's Derived Model theorem \cite{St09, StstattowfreeDMT}.} In fact, the aim is to build such determinacy models that reflect as much set theoretic strength of the universe as possible. This would, for example, be the case if $\HOD^M$ of our desired determinacy model $M$ has the same large cardinals as the universe $V$. Here $\HOD^M$ is the class of all hereditarily ordinal definable sets in $M$ and thus a canonical choice model inside any given model $M$ of determinacy. More precisely, a typical and well-known conjecture that one could approach via the core model induction technique is the following.  

\begin{conjecture}
    Suppose the Proper Forcing Axiom ($\PFA$) holds, let $\kappa \geq \omega_2$, and let $g$ be $\Col(\omega, \kappa)$-generic over $V$. Then there is a superstrong cardinal in $\HOD^{L(\uB_g, \bR_g)}$.
\end{conjecture}

Currently, this conjecture is open even for a Woodin limit of Woodin cardinals instead of a superstrong cardinal. To outline why $\Sealing$ is an obstruction in attempts to prove this conjecture we describe the standard approach to obtain large cardinals in $\HOD^{L(\uB_g, \bR_g)}$ via the core model induction. The key is the failure of a covering principle for $\HOD^{L(\uB_g, \bR_g)}$. For $\Theta$ the least ordinal that is not a surjective image of the reals, consider $\mathcal{H}^- = \HOD^{L(\uB_g, \bR_g)}|\Theta$ and let $\mathcal{H}$ be a canonical\footnote{We omit the technical definition of $\mathcal{H}$ here for simplicity. The interested reader can consult the introduction of \cite{SaTrConSealing} for a more detailed outline.} one-cardinal extension of $\mathcal{H}^-$. A typical attempt to prove that $\PFA$ implies that there is an inner model with a Woodin limit of Woodin cardinals goes as follows: We do core model induction below a cardinal $\kappa$ that is a measurable limit of strong cardinals (a singular strong limit cardinal $\kappa$ in fact suffices), so let $g$ be $\Col(\omega, {<}\kappa)$-generic over $V$. Assume that there is no Woodin limit of Woodin cardinals in $\mathcal{H}$. Then, the usual procedure in a core model induction argument, would be to show that $\mathsf{UB-Covering}$ fails and use this to show that there is a universally Baire set which is not in $\uB_g$, an obvious contradiction. Here $\mathsf{UB-Covering}$ is the statement that $\cf^V(\Ord \cap \, \mathcal{H}) \geq \kappa$. But, surprisingly, Sargsyan and Trang showed in \cite{SaTrConSealing} (see Theorem \ref{thm:ConSealing}) that this setup does not work anymore. It is not true that either there is no Woodin limit of Woodin cardinals in $\mathcal{H}$ or $\mathsf{UB-Covering}$ must hold.


\section{Conclusion}

So which scenario is true? Does $\Sealing$ hold in the universe? Does Woodin's axiom $\UltimateL$ hold? These are more general structural questions behind Cantor's Continuum Hypothesis and Gödel's Program that we need to answer. Answering these question would not only solve local problems such as whether there is a set of reals of size between the size of $\mathbb{N}$ and the size of $\bR$ but they give a much broader perspective of the mathematical universe. For example, if $\UltimateL$ holds, it not only implies $\CH$ but it gives an explanation why $\CH$ holds as it implies that there is a canonical well-ordering of the reals.

On the other hand the question whether $\Sealing$ or $\UltimateL$ hold lead us directly to the boundaries of our current knowledge and maybe even to the boundaries of mathematics. They are connected to even more fundamental questions such as: What is mathematical truth? And this is a question that leads us, as shown by Gödel in his famous incompleteness theorems, beyond the framework of mathematical proofs.

\bibliographystyle{plain}
\bibliography{References}

\end{document}